\documentclass[11pt]{article} 
\usepackage{amssymb}
\usepackage{latexsym,amsmath,amsthm,amsfonts}

\newtheorem{prop}{Proposition}
\newtheorem{lem}{Lemma}
\newtheorem{thm}{Theorem}
\newtheorem{cor}{Corollary}

\newcommand{\be}{\begin{equation}}
\newcommand{\ee}{\end{equation}}

\newcommand{\deq}{\stackrel{\scriptscriptstyle\triangle}{=}}

\newcommand{\B}{\mbox{$\cal B$}}
\newcommand{\C}{\mathcal{C}}
\newcommand{\sigS}{\mbox{$\cal S$}}
\newcommand{\sigB}{\mbox{$\cal B$}}
\newcommand{\D}{\mbox{$\cal D$}}

\newcommand{\X}{\mbox{$\cal X$}}

\begin{document}

\bibliographystyle{plain}

\title{Uniform Approximation and Bracketing Properties \\ of VC classes}

\author{Terrence M.\ Adams
\thanks{Terrence Adams is with the Department of Defense, 
9800 Savage Rd. Suite 6513, Ft. Meade, MD 20755} \ 
and Andrew B.\ Nobel 
\thanks{Andrew Nobel is with the Department of Statistics and 
Operations Research, University of North Carolina, Chapel Hill,
NC 27599-3260.  Email: nobel@email.unc.edu}}

\date{July 2010}

\maketitle

\begin{abstract}

We show that the sets in a family with finite VC
dimension can be uniformly approximated within
a given error by a finite partition.  
Immediate corollaries include 
the fact that VC classes have finite bracketing numbers, satisfy
uniform laws of averages under strong dependence, and exhibit
uniform mixing.  Our results are based on recent work concerning
uniform laws of averages for VC classes under ergodic sampling.
\end{abstract}

\newpage

\section{Introduction}

Let $\X$ be a complete separable metric space with Borel
sigma field $\sigS$, and let $\C \subseteq \sigS$ be a family 
of measurable sets.   For each finite set $D \subseteq \X$,
let $\{ C \cap D : C \in \C \}$ be the collection of subsets of
$D$ induced by the members of $\C$.
The family $\C$ is said to be a 
Vapnik-Chervonenkis (VC) class if there is a finite integer
$k$ such that
\be
\label{finvc}
| \{ C \cap D : C \in \C \} | \, < \, 2^k \  
\mbox{for every $D \subseteq \X$ with $|D| = k$} .
\ee
Here and in what follows $| \cdot |$ denotes cardinality.
The smallest $k$ for which (\ref{finvc}) holds is known
as the VC-dimension of $\C$.  
Classes of sets having finite VC-dimension play a central role in the 
theory of machine learning and empirical processes ({\it c.f.} 
\cite{Poll84, VaaWel96, DevGyoLug96, Dudley99}).

\subsection{Principal Result}

Let $\mu$ be a probability measure on $(\X,\sigS)$, and 
let $\pi$ be a finite, measurable partition
of $\X$.  For every set $C \in \C$,
the $\pi$-boundary of $C$, denoted $\partial(C : \pi)$, is the union of all the
cells in $\pi$ that intersect both $C$ and its complement
with positive probability.  Formally,
\[
\partial(C : \pi)
\ = \ 
\cup \, \{ A \in \pi : \mu(A \cap C) > 0 \mbox{ and }  \mu(A \cap C) > 0 \} .
\]
Note that $\partial(C : \pi)$ depends on $\mu$; this dependence
is suppressed in our notation.
Of interest here is the existence of a fixed finite partition $\pi$ such that the
measure of the boundary $\partial(C : \pi)$ is small for every
set $C$ in $\C$.  In general, the
existence of a uniformly approximating partition depends on the
family $\C$ and the measure $\mu$.  
Our main result shows that VC classes possess this uniform
approximation property, regardless of the measure $\mu$.

\vskip.2in

\begin{thm}
\label{BdyThm}
Let $\mu$ be a probability measure on $(\X, \sigS)$.
If $\C$ is a VC-class, then for every $\epsilon > 0$ there 
exists a finite measurable partition $\pi$ of $\X$ such that
\be
\label{approx}
\sup_{C \in {\cal C} } \mu(\partial(C : \pi) ) \, < \, \epsilon .
\ee
\end{thm}

\vskip.2in

Several corollaries of Theorem \ref{BdyThm} are discussed in the
next section.  The proof of Theorem \ref{BdyThm} is presented in 
Section \ref{PfBdyThm}.

\section{Corollaries of Theorem \ref{BdyThm}}

Here we present several immediate
corollaries of Theorem \ref{BdyThm} that may be of
independent interest.

\subsection{Bracketing of VC Classes}

Let $\mu$ be a probability measure on $(\X, \sigS)$.
For each pair of sets $A, B \in \sigS$, the bracket $[A,B]$ consists of all
those sets $C \subseteq \X$ such that $A \subseteq C \subseteq B$.  If $A$ 
is not a subset of $B$, then $[A,B]$ is empty.  The bracket
$[A,B]$ is said to be an $\epsilon$-bracket if $\mu(B \setminus A) \leq \epsilon$.
The bracketing number $N_{[ \, ]}(\epsilon,\C,\mu)$ 
of a family $\C \subseteq \sigS$
is the least number 
of $\epsilon$-brackets needed to cover $\C$.  Note that the 
sets defining the minimal brackets need not be elements of $\C$.  

\vskip.1in

\begin{cor}
\label{Bracket}
Let $\mu$ be any probability measure on $(\X, \sigS)$.
If $\C$ is a countable VC-class, then $N_{[ \, ]}(\epsilon,\C,\mu)$ 
is finite for every $\epsilon > 0$.
\end{cor}

\vskip.1in

\noindent
{\bf Remark:} Using routine arguments, the assumption that $\C$
is countable can be replaced by the weaker assumption that there exists
a countable sub-family $\C_0 \subseteq \C$ such that
the indicator function of every set in $\C$ is the pointwise limit of 
the indicator functions of sets in $\C_0$.

\vskip.15in

\noindent
{\bf Proof:} 
Fix a probability measure $\mu$ and $\epsilon > 0$.  
Let $\pi = \{A_1,\ldots,A_m\}$ be a finite measurable 
partition of $\X$ such that (\ref{approx}) holds, and assume
without loss of generality that each set $A_j$ has positive
$\mu$-measure.  
Let $A_j$ be an element of $\pi$.  For each $C \in \C$, 
remove points in $C$ from $A_j$ if $\mu(A_j \cap C) = 0$,
and remove points in $C^c$ from $A_j$ if  
$\mu(A_j \cap C^c) = 0$.  Denote the resulting
set by $B_j$.  Clearly $B_j \subseteq A_j$ and, as $\C$
is countable, $\mu(A_j \setminus B_j) = 0$.  The definition
of $B_j$ ensures that for each $C \in \C$ exactly one of the
following relations holds: $B_j \subseteq C$, 
$B_j \subseteq C^c$, or $\mu(B_j \cap C) \cdot \mu(B_j \cap C^c) > 0$.
Let $B_0 = \X \setminus \cup_{j=1}^m B_j$, and define the partition
$\pi' = \{ B_0, B_1, \ldots, B_m \}$.  Given $C \in \C$ let
$C_l = \cup \{ B \in \pi' : B \subseteq C \}$ and  
$C_u = \cup \{ B \in \pi' : B \cap C \neq \emptyset \}$.
A straightforward argument shows that 
$C_l \subseteq C \subseteq C_u$, and that
$\mu(C_u \setminus C_l) = 
\mu(\partial(C : \pi')) = \mu(\partial(C : \pi)) < \epsilon$.  
It follows that $\Theta = \{ [C_l,C_u] : C \in \C\}$ is a collection of 
$\epsilon$-brackets covering $\C$.   The cardinality of $\Theta$
is at most $2^{2 | \pi |}$.

\subsection{Uniform Laws of Large Numbers}

Let $X_1, X_2, \ldots$ be a stationary ergodic 
process taking values in $(\X,\sigS)$ with
$X_i \sim \mu$.  The ergodic theorem ensures that, for
every $C \in \sigS$, the sample averages 
$n^{-1} \sum_{i=1}^n I_C(X_i)$ converge with probability
one to $\mu(C)$.  For VC classes and i.i.d.\ sequences
$\{X_i\}$ this convergence is known to be uniform over $\C$
\cite{VC71}.
Using Corollary \ref{Bracket} it is easy to show that this
uniform convergence extends to ergodic processes as well.

\begin{thm}
\label{Ergodic}
If $\C$ is a countable VC-class of sets and 
$X_1, X_2, \ldots \in \X$ is a stationary ergodic 
process with 
$X_i \sim \mu$, then
\[
\sup_{C \in {\cal C}} 
\left| \frac{1}{n} \sum_{i=1}^n I_C(X_i) - \mu(C) \right|
\ \to \ 0
\]
with probability one as $n$ tends to infinity.
\end{thm}

\noindent
{\bf Proof:} This follows easily from Corollary \ref{Bracket}
and the Blum DeHardt law of large numbers ({\it c.f.} \cite{VaaWel96}),
which establishes that families with finite bracketing numbers 
have the Glivenko Cantelli property.

\vskip.1in

The uniform strong law in
Theorem \ref{Ergodic} was established in \cite{AdNob09}
using arguments similar to those forTheorem \ref{BdyThm}.  
Analogous uniform strong laws for VC major and VC graph 
classes are given in \cite{AdNob09}, while \cite{AdNob10}
contains uniform strong laws for classes of functions having
finite gap (fat shattering) dimension.  See these papers 
for a discussion of earlier and related work.

\subsection{Uniform Mixing Conditions in Ergodic Theory}

Let $T$ be an ergodic $\mu$-measure preserving transformation 
of $(\X,\sigS)$.  $T$ is said to be strongly mixing if for each pair 
$A$, $B$ of measurable sets, 
$\lim_{n \to \infty} \mu (A\cap T^{-n}B) = \mu (A)\mu (B).$ 
Theorem \ref{BdyThm} can be applied to show that strong mixing occurs uniformly over a countable VC class. 

\begin{prop}
\label{usm}
If $\C \subseteq \sigS$ is a countable VC-class of measurable sets, and 
$T$ is a strongly mixing transformation, then
\[
\lim_{n\to \infty} 
\sup_{A,B \in {\cal C}} 
\left| \mu (A\cap T^{-n}B) - \mu (A)\mu (B) \right| = 0.
\]
\end{prop}

\noindent
{\bf Proof:} 
Given $\epsilon > 0$, let $\pi$ be a finite partition such that 
$\sup_{C\in \C} \mu (\partial (C:\pi)) < \epsilon$. 
Choose a natural number $N$ such that for $n \geq N$ and 
each pair $D_1, D_2 \in \pi$, 
\[
|\mu (D_1 \cap T^{-n} D_2) - \mu (D_1) \mu (D_2)| \, < \,
\epsilon \, \mu (D_1) \mu (D_2).
\] 
For every measurable set $A$ let
$\overline{A} = \cup \{ D \in \pi: \mu (D \cap A) > 0 \}$ and 
$\underline{A} = \cup \{ D \in \pi: D \subset A\}$ be,
respectively, upper
and lower approximations of $A$ derived from the cells of $\pi$. 
Note that if $A,B$ are measurable sets satisfying 
$\overline{A}=\underline{A}$ and $\overline{B}=\underline{B}$, 
then 
\begin{eqnarray*}
|\mu (A \cap T^{-n}B) - \mu(A) \mu(B)| 
& = & 
|\sum_{D \subseteq \overline{A}} \sum_{D' \subseteq \overline{B}} 
\mu (D \cap T^{-n} D') - 
\sum_{D \subseteq \overline{A}} \sum_{D' \subseteq \overline{B}} 
\mu (D) \mu (D')| \\
& \leq &
\sum_{D \subseteq \overline{A}} \sum_{D' \subseteq \overline{B}} 
|\mu (D \cap T^{-n} D') - \mu(D) \mu (D')| \\ 
& < & 
\sum_{D \subseteq \overline{A}} \sum_{D' \subseteq \overline{B}} 
\epsilon \, \mu(D) \mu(D') 
\ \leq \ \epsilon \mu(A) \mu(B)
\ \leq \ \epsilon.
\end{eqnarray*}
Suppose now that $A,B$ are sets in $\C$. 
Then for $n\geq N$, 
\begin{eqnarray*}
\lefteqn{ |\mu (A\cap T^{-n}B) - \mu (A)\mu (B)|} \\ 
& = & 
|\mu (A\cap T^{-n}B) 
\, \pm \, \mu (A\cap T^{-n}\overline{B}) 
\, \pm \, \mu (\overline{A}\cap T^{-n}\overline{B})  
\, \pm \, \mu (\overline{A})\mu (\overline{B}) 
\, \pm \mu \, (\overline{A}) \mu (B) 
\, - \, \mu (A)\mu (B)| \\
& \leq &
2 \mu (\overline{B}\setminus B) \, + \, 2 \mu (\overline{A}\setminus A) 
\, + \, 
|\mu (\overline{A}\cap T^{-n}\overline{B}) - 
\mu (\overline{A})\mu (\overline{B})| \\ 
& < & 
5\epsilon ,
\end{eqnarray*}
where the first inequality follows from the triangle inequality,
and the second follows from the previous two displays.
As $A, B \in \C$ and $\epsilon >0$ were arbitrary, 
Theorem \ref{usm} follows.

\vskip.1in

A similar argument can be used to show that 
any weak mixing transformation satisfies 
uniform convergence over countable VC classes. 
A measure preserving transformation $T$ is weak mixing 
if given measurable sets $A$ and $B$, 
$$\lim_{n\to \infty} \frac{1}{n} \sum_{i=0}^{n-1} 
|\mu (A\cap T^{-i}B) - \mu (A)\mu (B)| = 0.$$ 

\begin{prop}
\label{uwm}
If $\C$ is a countable VC-class of measurable sets and 
$T$ is a weakly mixing transformation, then
\[
\lim_{n\to \infty} 
\sup_{A,B \in {\cal C}} 
\frac{1}{n} \sum_{i=0}^{n-1} 
|\mu (A\cap T^{-i}B) - \mu (A)\mu (B)| = 0. 
\]
\end{prop}

\section{Proof of Theorem \ref{BdyThm}}
\label{PfBdyThm}

The proof of Theorem \ref{BdyThm} follows arguments used 
in \cite{AdNob09} to establish uniform laws of large numbers
for VC classes under ergodic sampling, and we make use of 
several auxiliary results from that paper in what follows.

\subsection{Joins and the VC dimension}

\noindent
{\bf Definition:} The join of $k$ sets 
$A_1,\ldots, A_k \subseteq [0,1]$, 
denoted $J = \bigvee_{i=1}^k A_i$, is the partition consisting of all 
{\em non-empty} intersections 
$\tilde{A}_1 \cap \cdots \cap \tilde{A}_k$ 
where $\tilde{A}_i \in \{ A_i, A_i^c \}$ for $i = 1,\ldots,k$.  

\vskip.2in

Note that $J$ is a finite partition of $[0,1]$. 
The join of  $A_1,\ldots, A_k$ is said to be full if it has
(maximal) cardinality $2^k$. 
The next Lemma (see \cite{Mat02, AdNob09}) 
makes an elementary connection between full joins
and the VC dimension.  

\vskip.1in

\begin{lem}
\label{Join}
Let $\C$ be any collection of subsets of $\X$.  If for some $k \geq 1$ 
there exists a collection $\C_0 \subseteq \C$ of $2^k$ sets having 
a full join, then VC-dim$(\C) \geq k$.
\end{lem}

\vskip.1in

The proof given here establishes that the approximating partition 
$\pi$ is measurable $\sigma(\C)$.  A simple counterexample shows that 
it is not sufficient for the elements of $\pi$ to belong to 
$\bigcup_{n=1}^{\infty} \sigma(C_1, C_2, \ldots, C_n)$. 
To see this, let $\X=[0,1]$ and let $\lambda$ be Lebesgue measure. 
Let $a_1, a_2, \ldots > 0$ be a sequence of numbers such that 
$s=\sum_{n=1}^{\infty} a_n < 1$. 
Let $s_n = \sum_{i=1}^n a_i$ for $n \geq 1$ and let $s_0=0$. 
Define $C_n = [s_{n-1},s_n)$ for $n\geq 1$. 
Clearly, the VC-dimension of the class $\{C_1,C_2,\ldots\}$ 
equals 1, since its constituent sets are disjoint.
Define $J_n = C_1 \vee C_2 \vee \ldots \vee C_n$. 
Then $A_n=[s_n,1]$ is a
single element in $J_n$ with measure 
$1-s_n > 1-s > 0$.  Moreover, both $A_n\cap C_{n+1}$ 
and $A_n\cap C^{\prime}_{n+1}$ have positive measure,
so that
$\mu (\partial(C_{n+1} : J_n)) > 1-s$ for $n \geq 1$.

\subsection{Reduction to the Unit Interval} 

Fix a probability measure $\mu$ on $(\X,\sigS)$ and 
let $\C \subseteq \sigS$ have finite VC dimension.
It follows from standard results on the $L_p$-covering
numbers of VC classes ({\em c.f.} Theorem 2.6.4 of 
\cite{VaaWel96}) that there exists a countable sub-family
$\C_0$ of $\C$ such that
\[
\inf_{C' \in {\cal C}_0} \mu(C' \triangle C) = 0
\]
for each $C \in \C$.
An elementary argument then shows that, for every finite partition $\pi$,
\[
\sup_{C \in {\cal C}} \mu( \partial(C : \pi)) \ = \ 
\sup_{C \in {\cal C}_0} \mu( \partial(C : \pi)) ,
\]
and we may therefore assume that $\C$ is countable.
Let $\X_0 = \{ x : \mu(\{x\}) > 0 \}$ be the set of atoms
of $\mu$ and let 
$\mu_0(A) = \mu( A \cap {\cal X}_0 )$ be the
atomic component of $\mu$.  As $\X_0$ is countable,
it is easy to see that
\[
\inf_{\pi \in \Pi} \, \sup_{C \in {\cal C}} \mu_0( \partial(C : \pi))
\ = \ 0 ,
\]
and we may therefore assume that $\mu$ is non-atomic.

Following the proof in \cite{AdNob09}, we make two
further reductions.  
Let $\lambda(\cdot)$ be Lebesgue measure on the unit interval
$[0,1]$ equipped with its Borel subsets $\B$.
Using the existence of a measure-preserving
isomorphism between $(\X, \sigS, \mu)$ and 
$([0,1], \B, \lambda)$ ({\it c.f.} \cite{Royd88}) a straightforward
argument ensures that
we lose no generality in assuming that $\X = [0,1]$, $\mu = \lambda$, 
and that $\C \subseteq \B$ is a countable family with finite
VC dimension.  Using an additional isomorphism described 
in Lemma 6 of \cite{AdNob09} we may further assume that each 
element of $\C$ is a finite union of intervals.

Based on the reductions above,
Theorem \ref{BdyThm} is a corollary of the following result.

\begin{thm}
\label{Bdy-UI}
Let $\C \subseteq \B$ be a countable VC class, each of
whose elements is a finite union of intervals.
For every $\epsilon > 0$ there exists a finite
partition of $[0,1]$ such that
\[
\sup_{C \in {\cal C} } \lambda(\partial(C : \pi) ) \, < \, \epsilon .
\]
\end{thm}

\noindent
{\bf Remark:}  The proof of Theorem \ref{Bdy-UI} follows
the proof of Proposition 3 from
\cite{AdNob09}.
Beginning with the assumption that the conclusion of
the theorem is false, we construct, in a step-wise fashion,
a sequence of ``splitting sets'' $R_1, R_2, \ldots \subseteq [0,1]$
from the sets in $\C$.  
At the $k$th stage the splitting set $R_k$ is obtained from 
a sequential procedure that makes use of the splitting sets 
$R_1,\ldots, R_{k-1}$ produced at previous stages.  The splitting
sets are then used to identify finite, but arbitrarily large, collections of sets in 
$\C$ having full join.   The existence of these collections implies
that $\C$ has infinite VC dimension by Lemma \ref{Join}. 

\vskip.1in

\noindent 
{\bf Proof of Theorem \ref{Bdy-UI}:} 
Suppose to the contrary that there exists an $\eta > 0$
such that
\be
\label{etaineq}
\sup_{C \in {\cal C}} \lambda(\partial(C : \pi) ) 
 \ > \eta \, \mbox{ for every } \, \pi \in \Pi.
\ee
For $n \geq 1$ let $\D_n = \{ [\, k \, 2^{-n}, (k+1) \, 2^{-n}] : 0 \leq k \leq 2^n-1 \}$
be the set of closed dyadic intervals of order $n$.  

\vskip.1in

\noindent
{\bf Stage 1.}
Let $C_1(1)$ be any set in $\C$.  Suppose 
that sets $C_1(1),\ldots,C_1(n) \in \C$ have already been selected, and let  
$J_1(n) = \D_n \vee C_1(1) \vee \cdots \vee C_1(n)$.
It follows from (\ref{etaineq}) that there is a set $C_1(n+1) \in \C$
such that $G_1(n) = \partial(C_1(n+1) : J_1(n))$ has measure greater
than $\eta$.  
Let $J_1(n+1) = \D_{n+1} \vee C_1 \vee \cdots \vee C_{n+1}$ and
continue in the same fashion.  The sets $\{ G_1(n) \}$
are naturally associated with a tight family of sub-probability
measures $\{ \lambda_n (\cdot) = \lambda (\cdot \cap G_1(n)) \}$.   
There is therefore a subsequence
$\{ \lambda_{n_1(r)} \}$ that converges weakly to a sub-probability $\nu_1$ on
$([0,1],\sigB)$.  It is easy to see that $\nu_1$ is absolutely continuous 
with respect to $\lambda$ and that
\[
\nu_1([0,1]) 
\ \geq \ 
\limsup_{r \to \infty} \lambda_{n_r} ([0,1]) 
\ \geq \ 
\eta .
\]
The Radon-Nikodym derivative $d\nu_1 / d\lambda$ is well defined, 
and is bounded above by 1.  
Define the splitting set $R_1 = \{ x: (d\nu_1 / d\lambda)(x) > \eta/2 \}$.  From the previous 
remarks it follows that 
\begin{eqnarray}
\eta
& \leq &
\nu_1([0,1]) 
\ = \
\int_0^1 \frac{d\nu_1}{d\lambda} \, d\lambda 
\label{Rlbd}
\ \leq \ 
\int_{R_1} 1 d\lambda + \int_{R_1^c} \eta/2 \, d\lambda
\ \leq \ 
\lambda(R_1) + \eta/2 ,
\end{eqnarray}
and therefore $\lambda(R_1) \geq \eta/2$.  

\vskip.2in

\noindent
{\bf Subsequent stages.}
In order to construct the splitting set $R_k$ at stage $k$,  
let $C_k(1)$ be any element of $\C$,
and suppose that $C_k(2), \ldots, C_k(n)$ have already been
selected.  Define the join
\be
\label{joindef}
J_k(n) \ = \ \D_n \vee \bigvee_{j=1}^{k-1} R_j \vee \bigvee_{i=1}^n C_k(i) .
\ee
By (\ref{etaineq}) there exists a set $C_k(n+1) \in \C$ such that 
$G_k(n) = \partial(C_k(n+1 : J_k(n))$ has measure greater
than $\eta$.   This process continues as in stage 1.  
As before, there is a sequence of integers
$n_k(1) < n_k(2) < \cdots$ such that the measures
$\lambda(B \cap G_k(n_k(r)))$ converge weakly 
to a sub-probability measure $\nu_k$ on $([0,1], \sigB)$ that is absolutely
continuous with respect to $\lambda(\cdot)$.  
Define $R_k = \{ x: (d\nu_k / d\lambda)(x) > \delta \}$.  

\vskip.2in

\noindent
{\bf Construction of Full Joins.}
Fix an integer $L \geq 2$.   As the measures of the sets $R_k$ 
are bounded away from 
zero, there exist positive integers $k_1 < k_2 < \ldots < k_{L}$ such that                                                                                                                                                               
$\lambda(\bigcap_{j=1}^{L} R_{k_j}) > 0$.   Suppose without loss of 
generality that $k_j = j$, and define the intersections 
\[
Q_r = \bigcap_{j=1}^{L-r} R_j
\] 
for $r = 0, 1, \ldots, L-1$.  
Note that $Q_0 \subseteq Q_1 \subseteq \cdots \subseteq Q_{L-1}$.
We show that
there exist sets $D_1, D_2, \ldots, D_{L-1} \in \C$ such that,
for $l = 1,\ldots,L-1$, 
\begin{enumerate}

\item[(i)] the join
$K_l = D_1 \vee D_2 \vee \cdots \vee D_l$ has cardinality
$|K_l| = 2^l$, and 

\item[(ii)] $B^o \cap Q_l$ is non-empty
for each $B \in K_l$, where $B^o$ denotes the interior of $B$.

\end{enumerate}
We proceed by induction, beginning with the case $l = 1$.   
Let $x_1$ be a Lebesgue point of $Q_0$, 
and let $\epsilon = \eta / 2(\eta + 2)$.  Then there exists  
$\alpha_1 > 0$ such that the interval 
$I_1 \deq (x_1 - \alpha_1, x_1 + \alpha_1)$ satisfies
\be
\label{Int1prop}
\lambda(I_1 \cap Q_0) 
\ \geq \ 
(1 - \epsilon) \lambda(I_1) 
\ = \ 
2 \alpha_1 (1 - \epsilon) .
\ee
It follows from
the last display and the definition of $R_L \supseteq Q_0$ that
\be
\label{nukineq}
\nu_{L}(I_1 \cap R_{L})
\ = \ 
\int_{I_1 \cap R_L} \frac{ d\nu_{L} }{d\lambda} \, d\lambda
\ > \ 
\frac{eta}{2} \, \lambda(I_1 \cap R_L)
\ \geq \ 
\alpha_1 (1 - \epsilon) \eta .
\ee
Let $\{ n_L(r) : r \geq 1 \}$ be the subsequence used to define
the sub-probability $\nu_L$. 
As $I_1$ is an open set, the portmanteau theorem and (\ref{nukineq})
imply that
\[
\liminf_{r \to \infty} 
\lambda( I_1 \cap G_{L}(n_L(r)) )
\ \geq \ 
\nu_{L} (I_1)
\ \geq \ 
\nu_{L} (I_1 \cap R_L)
\ > \ 
\alpha_1 (1 - \epsilon) \eta .
\]
Choose $r$ sufficiently large so that 
$\lambda( I_1 \cap G_{L}(n_L(r)) ) \, > \, 
\alpha_1 (1 - \epsilon) \eta$ and
$2^{ - n_L(r) } < \eta \, \alpha_1 / 8$.  We require
the following lemma from \cite{AdNob09}.

\vskip.16in
 
\begin{lem}
\label{Alem}
There exists a cell $A$ of $J_L(n_L(r))$ such that
$A \subseteq \partial(C_L(n_L(r) + 1): J_L(n_L(r))$,
$A \subseteq I_1$ and $\lambda(A \cap Q_1) > 0$.  Moreover,
$A$ is contained in $Q_1$. 
\end{lem}

\vskip.1in

Let $D_1 = C_L(n_L(r)+1) \in \C$, and let $A$ be
the set identified in Lemma \ref{Alem}.  By definition of 
the boundary, $\lambda(A \cap D_1) > 0$ and
$\lambda(A \cap D_1^c) > 0$ and therefore
$\lambda(Q_1 \cap D_1) > 0$ and
$\lambda(Q_1 \cap D_1^c) > 0$ as well.
As the Lebesgue measure of the boundary
$\overline{D}_1 \setminus D_1^o$ of $D_1$ is zero, 
assertion (ii) above follows.

Suppose now that we have identified sets 
$D_1,\ldots,D_l \in \C$, with $l \leq L-2$, 
such that (i) and (ii) hold.
Let the join $K_l = \{ B_j : 1 \leq j \leq 2^l \}$, and for each $j$
let $x_j \in B_j^o \cap Q_l$.
Select $\alpha_{l+1} > 0$ such that for each $j$ the interval
$I_j \deq (x_j - \alpha_{l+1}, x_j + \alpha_{l+1})$
is contained in $B_j^o$ and satisfies
\[
\lambda(I_j \cap Q_l) 
\ \geq \ 
(1 - \epsilon) \lambda(I_j) 
\ = \
2 \alpha_{l+1} (1 - \epsilon) . 
\]
To simplify notation, let $\kappa = L - l$.
Let $\{ n_{\kappa} (r) : r \geq 1 \}$ 
be the subsequence used to define the sub-probability $\nu_{\kappa}$. 
For each interval $I_j$, 
\[
\liminf_{r \to \infty} 
\lambda( I_j \cap G_{\kappa}(n_{\kappa} (r) ) )
\ \geq \ 
\nu_{\kappa} (I_j)
\ \geq \ 
\nu_{\kappa} (I_j \cap R_{\kappa})
\ > \ 
\alpha_{l+1} (1 - \epsilon) \eta ,
\]
where the last inequality follows from the previous display,
and the fact that $Q_l \subseteq R_{\kappa}$. 
Choose $r$ sufficiently large so that 
$\lambda( I_j \cap G_{\kappa}(n_{\kappa} (r)) ) >
\alpha_{l+1} (1 - \epsilon) \eta$ for each $j$,
and $2^{ - n_{\kappa}(r) } < \eta \, \alpha_{l+1} / 8$.  

By applying the Lemma \ref{Alem} to each interval $I_j$,
one may establish the existence of 
sets $A_j \in \partial(C_{\kappa}(n_{\kappa}(r) + 1): J_{\kappa}(n_{\kappa}(r))$ such that 
$A_j \subseteq I_j \subseteq B_j^o$, $\lambda(A_j \cap Q_{l+1}) > 0$, and
$A_j \subseteq  Q_{l+1}$.
Let $D_{l+1} = C_{\kappa}(n_{\kappa} (r)+1) \in \C$.  
Arguments like those for the case $l = 1$ above
show that for each $j$ the intersections
$A_j \cap D_{l+1}^o$ and $A_j \cap (D_{l+1}^c)^o$ are non-empty,
and the inductive step is complete.
Given any two dyadic intervals, they are disjoint, intersect
at one point, or one contains
the other.   Therefore, among the sets $D_1,\ldots,D_{L-1}$, 
at most one can be a dyadic interval; the remainder
are contained in $\C$.

\vskip.3in

\noindent
{\bf\Large  Acknowledgements} \\
The authors would like to thank Ramon van Handel for pointing
out an oversight in the proof of Corollary \ref{Bracket}.
The work presented in this paper was supported in part
by NSF grant DMS-0907177.

\end{document}